\documentclass[amssymb,twoside,12pt,amscd,reqno]{amsart}

\usepackage{latexsym,amsfonts}
\usepackage{amsmath,amsthm}          
\usepackage{amssymb}                 
\usepackage{amscd}                   
\usepackage{euscript}                

\theoremstyle{plain}
\numberwithin{equation}{section}
\newtheorem{thm}{Theorem}[section]

\newtheorem{cor}[thm]{Corollary}

\newtheorem{lemma}[thm]{Lemma}

\newcommand{\bi}{\begin{itemize}}
\newcommand{\ei}{\end{itemize}}
\newcommand{\bp}{\begin{proof}}
\newcommand{\ep}{\end{proof}}

\def\CC{\mathbb{C}}
\def\GGG{\mathbb{G}}
\def\PP{\mathbb{P}}
\def\QQ{\mathbb{Q}}

\def\ZZ{\mathbb{Z}}

\def\ka{\kappa}
\def\la{\lambda}

\def\si{\sigma}
\def\th{\theta}

\def\ze{\zeta}

\def\Ga{\Gamma}

\def\Si{\Sigma}
\def\Th{\Theta}

\def\A{\mathcal{A}}
\def\AA{\tilde{\A}}

\def\D{\mathcal{D}}
\def\Lb{\mathcal{L}}
\def\E{\mathcal{E}}
\def\F{\mathcal{F}}

\def\G{\mathcal{G}}

\def\HH{\mathcal{H}}

\def\J{\mathcal{J}}

\def\K{\mathcal{K}}
\def\M{\mathcal{M}}
\def\MM{\overline{M}}
\def\N{\mathcal{N}}

\def\S{\mathcal{S}}

\def\S{\mathcal{S}}
\def\T{\mathcal{T}}

\def\V{\mathcal{V}}

\def\X{\mathcal{X}}

\def\Z{\mathcal{Z}}

\def\ad{\text{ad}}
\def\AA{\text{A}}
\def\Aut{\text{Aut }}
\def\codim{\text{codim }}
\def\coker{\text{coker }}
\def\deg{\text{deg}}
\def\det{\text{det}}
\def\dim{\text{dim }}

\def\EEnd{\mathcal{E}nd}
\def\Ext{\text{Ext}}
\def\EExt{\mathcal{E}xt}
\def\h{\text{h}}
\def\H{\text{H}}

\def\Hom{\text{Hom}}

\def\id{\text{id}}
\def\im{\text{im }}
\def\ker{\text{ker }}

\def\Pic{\text{Pic}}
\def\PGL{\text{PGL}}
\def\Proj{\text{Proj}}

\def\R{\text{R}}

\def\Sym{\text{Sym}}

\def\dra{\dashrightarrow}

\def\Lra{\Longrightarrow}
\def\ra{\rightarrow}

\def\bo{\boxtimes}

\def\eset{\emptyset}

\def\setm{\setminus}

\def\ps{\vspace{4pt}}

\begin{document}

\title{Rational Families of Vector Bundles on Curves, I}
\author[Castravet]{Ana-Maria Castravet}
\address{Institute for Advanced Study \\ Princeton, NJ 08540}
\email{noni@ias.edu}
\thanks{This material is based upon work partially supported by the National
Science Foundation under agreement No. DMS-9729992.}

\begin{abstract}
Let $C$ be a smooth complex projective curve of genus $g\geq2$ and let
$M$ be the moduli space of rank $2$, stable vector bundles on $C$, with fixed
determinant of degree $1$. For any $k\geq2$, we find two irreducible
components of the space of rational curves  of degree $k$ on $M$, both of
the expected dimension. One component, which we call the
\emph{nice component}, has the property that when $k$ is sufficiently
large, the general element of $\M$  is a very free curve. The second
component, that we call the \emph{almost nice component}, has the property
that the general element is a free curve
$f:\PP^1\ra M$ with $f^*T_M\cong O^g\oplus\A$, for some positive
vector bundle $\A$ on $\PP^1$. We prove that the maximal rationally
connected quotient of each component is the Jacobian $J(C)$ of the curve $C$.
\end{abstract}

\maketitle

\tableofcontents

\section{Introduction and background}

Let $C$ be a genus $g\geq2$ smooth projective curve over $\CC$.
Let $\xi$ be a degree $1$ line bundle on $C$ and let $M$ be the moduli
space of isomorphism classes of stable, rank $2$ vector bundles
$\E$ on C, with determinant $\det(\E)=\wedge^2(\E)\cong\xi$. Then $M$ is
a smooth, projective, irreducible variety of  dimension $(3g-3)$.
It is known that $\Pic(M)\cong\ZZ$ \cite{DN} and let $\Th$ be the ample
generator. In fact, $\Th$ is very ample \cite{BV}. By  a rational curve of
degree $k$ on $M$ we will always mean a non-constant morphism $f:\PP^1\ra M$,
such that the line bundle $f^*\Th$ on $\PP^1$ has degree $k$.
Our main goal in this paper is to study the Hilbert scheme
$\Hom_k(\PP^1,M)$ of morphisms  $f:\PP^1\ra M$ of degree $k$.

It is known that the canonical bundle
$K_M$ is $\Th^{-2}$ \cite{R}; hence, $M$ is a Fano variety. Moreover, $M$ is
a rational variety. There are some explicit descriptions of
$M$, that we will not use in what follows, but they are interesting to
mention. If $g=2$ then $M$ is isomorphic to a complete intersection of two
quadrics in $\PP^5$ \cite{N1}; more generally, for any $g\geq2$, if $C$ is
a hyperelliptic curve of genus $g$, then $M$ is isomorphic to the
Grassmanian of $(g-2)$-planes contained in a complete intersection of
two quadrics in $\PP^{2g+1}$ \cite{DR}.

\

Our motivation comes from the question of measuring the extent
to which the space of rational curves on a rationally connected variety is
itself rationally connected. Recall that a projective variety $X$ is
\emph{rationally connected}, if
there is a dense open $X^0\subset X$, such that for any two points $x_1$
and $x_2$ in $X^0$, there is a rational curve $\PP^1\ra X$ through $x_1$
and $x_2$. For curves and surfaces, this notion is equivalent to the
notion of rationality and unirationality; in higher dimensions,
rational connectivity is weaker. If $X$ is a smooth projective
variety, its
\emph{maximal rationally connected fibration (MRC)}
is a rational map $\psi:X\dra Z$ such that:
\bi
\item[i. ] the general fiber is rationally connected
\item[ii. ] for very general $z\in Z$, a rational curve in $X$ intersecting
the fiber $\psi^{-1}(z)$ is contained in $\psi^{-1}(z)$
\ei

The MRC fibration of a smooth projective variety is known to exist \cite{K}.
The pair $(Z,\psi)$ is unique up to a birational transformation.
The MRC fibration of a singular projective variety is defined to be 
the MRC fibration of some desingularization.

The question about the inductive behavior of rational connectivity for some
hypersurfaces of low degree has been the object of a series of
papers by Harris, Roth and Starr in \cite{HRS1}, \cite{HRS2}, \cite{HRS3}
\cite{HS}. Generally speaking, if $X$ is a smooth, projective, rationally
connected variety, we ask what are the irreducible components of the space
$\Hom_k(\PP^1, X)$, what are their dimensions and their MRC fibrations.
The answer for the case when $X$ is a projective space or a smooth quadric
of dimension at least three is given by Kim and Pandharipande in \cite{KP}.

\

In this paper (and in its follow-up \cite{C}), we answer these questions for
the case of moduli spaces $M$ of vector bundles. The problem of studying the
irreducible components of the space $\Hom_k(\PP^1, M)$, was considered before
in  \cite{Ki}, using a different method, for the case of $k=1$ and $k=2$.
Our methods are based on the classical ideas  of exhibiting a rank $2$ vector
bundle as the middle term of an extension of line bundles or an extension of
a skyscraper sheaf by a rank $2$ vector bundle.
One of the main results in this paper is the following:

\begin{thm}
For any integer $k\geq1$, there is a \emph{nice} irreducible component
$\M$ of the space $\Hom_k(\PP^1,M)$:
\bi
\item[i) ]it has the expected dimension $(2k+3g-3)$ and the general point is
unobstructed
\item[ii) ] the MRC  fibration of $\M$ is given by a rational map
$\M\dra J(C)$
\ei
\end{thm}

To prove the theorem for $k$ odd, we fix an integer $e\geq0$ and for every
degree $-e$ line bundle $\Lb$ on $C$, we consider exact sequences
\begin{equation}\label{ext_1}
0\ra\Lb\ra\E\ra\Lb^{-1}\otimes\xi\ra0
\end{equation}

If $e\geq0$ is an integer, for every degree $-e$ line bundle $\Lb$ on $C$,
denote with $V_{\Lb}$ the  vector space $\Ext^1(\Lb^{-1}\otimes\xi,\Lb)$,
parameterizing extensions of type (\ref{ext_1}). There is a rational map
$\PP(V_{\Lb})\dra M$ which associates to an extension (\ref{ext_1}) the
isomorphism class of the bundle $\E$.
We prove that it is defined outside a locus of codimension at least $2$. By
considering lines in the projective spaces $\PP(V_{\Lb})$, we obtain rational
curves of degree $(2e+1)$ on $M$. If we do this for all line bundles $\Lb$ of
degree $-e$, the rational curves obtained in this way fill up a whole
irreducible component of $\Hom_k(\PP^1,M)$, for $k=(2e+1)$.

To prove the theorem for $k$ even, we consider exact sequences
\begin{equation}\label{ext_2}
0\ra\E\ra\E'\ra O_D\ra0
\end{equation}
where $D\in\Sym^e(C)$ ($e\geq0$) and $\E$ is a rank $2$ stable vector
bundle, with $\det(\E)\cong\xi(-D)$. In a similar way as in the odd case,
if we denote with $V_{D,\E}$
the  vector space $\Ext^1(O_D,\E)$, parameterizing extensions of type
(\ref{ext_2}), we prove that there is a rational map $\PP(V_{D,\E})\dra M$,
that it is defined outside a locus of codimension at least $2$. By considering
lines in the projective spaces  $\PP(V_{D,\E})$, we obtain rational curves of
degree $2e$ on $M$. If we
do this for all $D\in\Sym^e(C)$ and all rank $2$ stable bundles $\E$ with
$\det(\E)\cong\xi(-D)$, the rational curves obtained in this way fill up a
whole irreducible component of $\Hom_k(\PP^1,M)$, for $k=2e$.

\ps

We ask whether the rational curves $f:\PP^1\ra M$ in the nice component are
free curves, that is if $f^*T_M$ is a non-negative vector bundle. The
following theorem states precisely when this is the case.

\begin{thm}\label{free_curves}
For any integer $k\geq1$, let $\M\subset\Hom_k(\PP^1,M)$ be the nice
component. We distinguish the following cases:
\bi
\item[i) ]If $k$ is odd, a general $f\in\M$ is a free curve
if and only if $k\geq(g-1)$
\item[ii) ]If $k$ is even, a general $f\in\M$ is a free curve
for any $k$
\ei
Moreover, if $k$ is sufficiently large, the general $f\in\M$ is a very free
curve.
\end{thm}

\ps

If we generalize the idea from the case when $k$ is odd, we obtain that,
at least in certain cases, there is at least another irreducible component.
In \cite{C}, we will give a different method, that helps to find and describe
all the irreducible components of the space $\Hom_k(\PP^1,M)$.

\begin{thm}
If $g$ is even, for integers $k\geq1$ that are divisible by $(g-1)$,
there is an \emph{almost nice} irreducible component $\M'$
(different from $\M$) of the space $\Hom_k(\PP^1,M)$:
\bi
\item[i) ]it has the expected dimension $(2k+3g-3)$
\item[ii) ] a general $f\in\M$ is a free curve $f:\PP^1\ra M$, such that
$f^*T_M\cong O^g\oplus\A$, where $\A$ is some positive vector bundle on
$\PP^1$
\item[iii) ] the MRC  fibration of $\M$ is given by a rational map
$\M\dra J(C)$
\ei
\end{thm}

\

We recall that if $X$ is a smooth projective variety and $f:\PP^1\ra X$ is a
morphism, there is a lower bound for the dimension of the irreducible
component of $\Hom(\PP^1, X)$, called the \emph{expected dimension} of
$\Hom(\PP^1, X)$ at the point $f$. It is given by the Euler characteristic
of the bundle $f^*T_X$ on $\PP^1$, which by Riemann-Roch is:
\begin{equation*}
\chi(\PP^1,f^*T_X)=\deg(f^*T_X)+\dim X=-f_*[\PP^1].K_X+\dim X
\end{equation*}

Moreover, the tangent space to $\Hom(\PP^1, X)$ at $f$ is isomorphic
to $\H^0(\PP^1, f^*T_X)$. A point $f$ is called \emph{unobstructed}
if $\H^1(\PP^1, f^*T_X)=0$. Note that it follows that $\Hom(\PP^1, X)$
is smooth at an unobstructed point $f$; hence, there is a unique
irreducible component of $\Hom(\PP^1, X)$, of the expected dimension,
that contains $f$.

\

The paper is organized as follows. In Section 2 we give the construction
about extensions of line bundles. In Section 3 we construct rational curves
in spaces of extensions of line bundles and we prove the existence of the
nice component for $k$ odd. In Section 4, we compute $f^*T_M$ for the
rational curves $f:\PP^1\ra M$ constructed in Section 3 and prove Theorem
\ref{free_curves}
for $k$ odd. We also prove the existence of the almost nice component.
In Section 5, we give the construction about extensions of skyscraper sheaves
by rank $2$ vector bundles. In Section 6 we construct rational curves in such
spaces of extensions and we prove the existence of the nice component for
$k$ even. In Section 7, we compute $f^*T_M$ for the rational curves
$f:\PP^1\ra M$ constructed in Section 6 and prove Theorem \ref{free_curves}
for $k$ even. In Section 8 we summarize our results for the case when $g=2$.

\ps

All our schemes, morphisms and products of schemes are over $\CC$, unless
otherwise stated. If $\E$ is a locally free sheaf on a scheme $X$, by
$\PP(\E)$ we mean $\Proj(\Sym(\E^*)$. When working with line bundles over
$\PP^1\times C$, if $\M$ is a line bundle on $\PP^1$ and $\N$ is a line bundle
on $C$, we use the notation $\M\bo\N$ for $p_1^*\M\otimes p_2^*\N$, where
$p_1$ and $p_2$ are the two projections.

\

\textbf{Acknowledgements.} I am very much indebted to Johan de Jong
for many suggestions and ideas. I also thank Steven
Kleiman, Mihnea Popa and Jason Starr, for helpful discussions.
Finally, I would like to thank my advisor, Joe Harris, for
introducing me to this beautiful subject and for the advice
and inspiration that he provided for me in the last couple of years.


\section{Extensions of line bundles}\label{ext_line_bundles}

Let $d$ be an integer and $\ze$ be a line bundle on $C$ of degree $d$.
Let $M_{\ze}$ be the coarse moduli scheme of semistable, rank $2$ vector
bundles on $C$, having determinant $\ze$. We construct maps from spaces of
extensions of line bundles on $C$ to $M_{\ze}$. The fact that we do this
for any degree $d$ line bundle $\ze$, and not only for $\xi$, will be
useful in the proof of Lemma \ref{unstable_locus_2_local}.

\ps

Let $e$ be an integer such that $e>-\frac{d}{2}$ and let $\Lb$ be a
line bundle on $C$ of degree $-e$. Consider extensions:
\begin{equation}\tag{$*$}
0\ra\Lb\ra\E\ra\Lb^{-1}\otimes\ze\ra 0.
\end{equation}
Then $\E$ is a rank $2$ vector bundle and such extensions are classified
by the vector space
\begin{equation*}
V_{\Lb}=\Ext^{1}_C(\Lb^{-1}\otimes\ze, \Lb)\cong\H^1(C,\Lb^2\otimes\ze^{-1})
\end{equation*}
Using Riemann-Roch, the vector space $V_{\Lb}$ has dimension
\begin{equation*}
\dim V_{\Lb}=-\deg(\Lb^2\otimes\ze^{-1})-1+g=2e+d+g-1
\end{equation*}

Clearly, any two nonzero elements $v,v'$ of $V_{\Lb}$ which differ by a
scalar define isomorphic vector bundles $\E$. Therefore the isomorphism
classes of non-trivial extensions as above are parametrized by the
projective space
$\PP(V_{\Lb})=\Proj(\Sym V_{\Lb}^*)$.

\ps

We now let $\Lb$ vary in the Picard variety $\Pic^{-e}(C)$ of line
bundles on $C$ of degree $-e$ and define a global parameter space for
extensions of type $(*)$. Let $\A$ be a {Poincar\'e} bundle on
$\Pic^{-e}(C)\times C$ and let $\pi_1,\pi_2$ be the two projections from
$\Pic^{-e}(C)\times C$. Define on $\Pic^{-e}(C)$ the relative extension sheaf
\begin{equation*}
\S=\EExt^1_{\Pic^{-e}(C)\times C|\Pic^{-e}(C)}
(\pi_1^*\A^{-1}\otimes\pi_2^*\ze,\pi_1^*\A)
\end{equation*}

Note that $\S$ is a locally free sheaf. Let $X$ be the projective bundle
$\PP(\S)$ and let
\begin{equation}\label{X}
p:X\ra\Pic^{-e}(C)
\end{equation}
be the corresponding map. (Note that $X$ depends on the choice of $\A$.)

Let $\nu_1,\nu_2$ be the two projections from $X\times C$. There is a
universal extension on $X\times C$:
\begin{equation}\label{univ_ext_1}
0\ra \nu_1^*O_X(1)\otimes p^*\A\ra\G\ra
p^*(\A^{-1})\otimes\nu_2^*\ze\ra0.
\end{equation}
It is universal in the sense that, if $x\in X$ and we let
$\Lb=p(x)\in\Pic^{-e}(C)$, when we restrict (\ref{univ_ext_1}) to
$\{x\}\times C$, we get an exact sequence:
\begin{equation}\label{ext_at_x}
0\ra\Lb\ra\G_x\ra\Lb^{-1}\otimes\ze\ra0
\end{equation}
whose class in $\PP(V_{\Lb})\cong p^{-1}(\{\Lb\})$ is $x$.


\subsection{The locus of unstable extensions}

Consider the space of extensions $X$. Let $Z\subset X$ (respectively
$Z'\subset X$) be the locus of those $x\in X$ for which the bundle $\G_x$
of (\ref{ext_at_x}) is unstable (respectively not stable). Since being
semistable (respectively stable) is an open condition, it follows that
$Z$ and $Z'$ are closed subsets of $X$. The following lemma gives an
estimate of their codimension.

\begin{lemma}\label{unstable_locus_1_local}
For each $\Lb\in\Pic^{-e}(C)$, there are projective irreducible varieties
$Z_{\Lb}$ and $Z'_{\Lb}$, such that
$Z_{\Lb}\subset Z'_{\Lb}\subset\PP(V_{\Lb})$, corresponding to the
unstable, respectively not stable, extensions in
$\PP(V_{\Lb})$ and we distinguish the following cases:
\bi
\item[i) ]If $d$ is odd, then $Z_{\Lb}=Z'_{\Lb}$ has codimension at least $g$.
\item[ii) ]If $d$ is even, then $Z_{\Lb}$ has codimension at least $(g+1)$ and
$Z'_{\Lb}$ has codimension at least $(g-1)$.
\ei
In both cases, when $e=1-\lceil\frac{d}{2}\rceil$, then $Z_{\Lb}=\eset$.
\end{lemma}

\bp
We analyze first the not stable locus $Z'_{\Lb}$.
The bundle $\E$ in ($*$) is not stable if and only if
there exists a line bundle $\Lb'$ on $C$ of degree $\lceil\frac{d}{2}\rceil$
and a non-zero morphism $\Lb'\ra\E$. Then the morphism
$\Lb'\ra\Lb^{-1}\otimes\ze$ is non-zero as well. This is because there
is no non-zero morphism $\Lb'\ra\Lb$ as
\begin{equation*}
\deg(\Lb')=\lceil
\frac{d}{2}\rceil>-e=\deg(\Lb)
\end{equation*}
Therefore, it follows that there is some effective divisor $D$ on $C$ of
degree $e+d-\lceil\frac{d}{2}\rceil$ such that
$\Lb'\cong\Lb^{-1}\otimes\ze(-D)$.
Let $\E'$ be the kernel of the composition morphism
\begin{equation*}
\E\ra\Lb^{-1}\otimes\ze\ra\Lb^{-1}\otimes\ze_{|D}
\end{equation*}
There is a commutative diagram with the two horizontal sequences exact:
\begin{equation*}
\begin{CD}
0@>>>\E'@>>>\E@>>>\Lb^{-1}\otimes\ze_{|D}@>>>0\\
@VVV @VVV @VVV @|@VVV\\
0@>>>\Lb^{-1}\otimes\ze(-D)@>>>\Lb^{-1}\otimes\ze@>>>
\Lb^{-1}\otimes\ze_{|D}@>>>0
\end{CD}
\end{equation*}
Using the snake lemma, we get that there is an exact sequence
\begin{equation}\label{split_seq}
0\ra\Lb\ra\E'\ra\Lb^{-1}\otimes\ze(-D)\ra0.
\end{equation}
Moreover, the following composition of morphisms is zero:
\begin{equation*}
\Lb^{-1}\otimes\ze(-D)\cong\Lb'\ra\E\ra\Lb^{-1}\otimes\ze_{|D}
\end{equation*}
Then $\Lb^{-1}\otimes\ze(-D)$ must map to the subbundle $\E'$ of $\E$ and
the exact sequence (\ref{split_seq}) is split.
We conclude that the vector $v\in V_{\Lb}=\Ext^1_C(\Lb^{-1}\otimes\ze,\Lb)$,
corresponding to an unstable extension, is in the kernel
$\Ext^1_C(\Lb^{-1}\otimes\ze\otimes O_D,\Lb)$ of the surjective map
\begin{equation*}
\Ext^1_C(\Lb^{-1}\otimes\ze,\Lb)\ra
\Ext^1_C(\Lb^{-1}(-D)\otimes\ze,\Lb)
\end{equation*}
Using duality, we have
\begin{equation*}
\Ext^1_C(\Lb^{-1}\otimes\ze\otimes O_D,\Lb)\cong
\H^0(C,\Lb^{-2}\otimes\ze\otimes K_C\otimes O_D)
\end{equation*}
Therefore, the non-stable extensions in $V_{\Lb}$ form a
$l:=\deg(D)=e+d-\lceil\frac{d}{2}\rceil$ dimensional linear subspace,
for each $D$ effective divisor of degree $l$. If we let $D$ vary in
$\Sym^l(C)$ we get that the locus of not stable extensions in $V_{\Lb}$
is at most $2l$-dimensional.

\ps

To make this precise, let $\D\subset\Sym^l(C)\times C$ be the
universal divisor and $p_1$ and $p_2$ be the two projections
from $\Sym^l(C)\times C$. Define the relative extension sheaf
\begin{equation*}
\F=\EExt^1_{\Sym^l(C)\times C|\Sym^l(C)}
(p_2^*(\Lb^{-1}\otimes\ze)\otimes O_{\D},p_2^*\Lb)
\end{equation*}

Note that $\F$ is a locally free sheaf; moreover, $\F$ is a subbundle
of the trivial bundle $O\otimes V_{\Lb}$ on $\Sym^l(C)$.
We have that the projective bundle $\PP(\F)$ is a closed
subvariety in $\Sym^l(C)\times\PP(V_{\Lb})$  and the not stable locus
$Z'_{\Lb}$ in $\PP(V_{\Lb})$ is given by the image of $\PP(\F)$ in
$\PP(V_{\Lb})$.
It follows that
\begin{equation*}
\codim_{\PP(V_{\Lb})}Z'_{\Lb}\geq
(2e+d+g-1)-2l=2\lceil\frac{d}{2}\rceil-d+g-1
\end{equation*}

If $d$ is odd, then $Z'_{\Lb}=Z_{\Lb}$ and its codimension is at least $g$.
If $d$ is even, the codimension of $Z'_{\Lb}$
is at least $g-1$. As for the unstable locus $Z_{\Lb}$, by a similar
argument, the bundle $\E$ is not semistable if and only if there exists
on $C$ a degree $\frac{d}{2}+1$ line bundle $\Lb'$ and a non-zero morphism
$\Lb'\ra\E.$ By degree considerations, the morphism
$\Lb'\ra\Lb^{-1}\otimes\ze$ is non-zero and we have
$\Lb'\cong\Lb^{-1}\otimes\ze(-D)$ for some effective divisor
$D$ on $C$ of degree $l:=e+\frac{d}{2}-1.$ The locus of unstable extensions
in $V_{\Lb}$ is a linear subspace of dimension at most $2l$. Hence, the
codimension of $Z_{\Lb}$ is at least $(2e+d+g-1)-2l=g+1$.

Note that, in both cases, if $l=0$
($e=1-\lceil\frac{d}{2}\rceil$), we have $Z_{\Lb}=\eset$.
\ep

\begin{cor}\label{unstable_locus_1}
There are projective irreducible varieties $Z$ and $Z'$, such that
$Z\subset Z'\subset X$, corresponding to the unstable, respectively
not stable, extensions in $X$ and we distinguish the following cases:
\bi
\item[i) ]If $d$ is odd, then $Z=Z'$ has codimension at least $g$ in $X$
\item[ii) ]If $d$ is even, then $Z$ (resp. $Z'$) has codimension at least
$(g+1)$ (resp. $(g-1)$) in $X$
\ei
In both cases, if $e=1-\lceil\frac{d}{2}\rceil$, then $Z=\eset$.
\end{cor}

It follows from Corollary \ref{unstable_locus_1} and from the coarse
moduli property of the moduli scheme $M_{\ze}$, that the universal bundle
$\G$ of (\ref{univ_ext_1}), restricted to $(X\setm Z)\times C$, induces
a well-defined morphism
\begin{equation}\label{ka}
\ka:X\setm Z\ra M_{\ze}
\end{equation}
which sends an extension $(*)$ to the class of the bundle $\E$ in $ M_{\ze}$.


\subsection{The global space of extensions of line bundles}

Let now $\ze$ vary in $\Pic^{d}(C)$ and denote
$V_{\Lb,\ze}=\Ext^1_C(\Lb^{-1}\otimes\ze,\Lb)$. In a similar way as in the
construction of (\ref{X}), there is a projective bundle
\begin{equation*}
q:\X\ra\Pic^{-e}(C)\times\Pic^{d}(C)
\end{equation*}
whose fiber at a point $(\Lb,\ze)\in\Pic^{-e}(C)\times\Pic^{d}(C)$
is canonically isomorphic to $\PP(V_{\Lb,\ze})$.

Let $\rho_1, \rho_2$ be the two projections from $\X\times C$. There is a
universal extension on $\X\times C$:
\begin{equation}\label{univ_ext_global}
0\ra q^*\A_1\otimes \rho_1^*O_{\X}(1)\ra\HH
\ra q^*({\A_1}^{-1}\otimes\A_2)\ra0
\end{equation}

As in Corollary \ref{unstable_locus_1}, there are projective irreducible
varieties $\Z$ and  $\Z'$, corresponding to the unstable, respectively not
stable, extensions in $\X$ and we have  $\Z\subset\Z'\subset\X$.

Let $M_d$ be the coarse moduli scheme of rank $2$, degree $d$, semistable
vector bundles of degree $d$. Then the universal bundle $\HH$ of
(\ref{univ_ext_global}) induces a morphism $\tilde{\ka}$, that is a global
version of the morphism $\ka$ of (\ref{ka}):
\begin{equation}\label{ka_global}
\tilde{\ka}:\X\setm\Z\ra M_d
\end{equation}


\

\textbf{Notation.}
From now on, we let $\xi$ to be a degree $1$ line bundle on $C$ and
we denote by $M$ the moduli space of stable, rank $2$ vector bundles on $C$,
with determinant $\xi$.


\subsection{The morphism $\ka$}

Consider the morphism (\ref{ka}) for the case when $\ze=\xi$:
\begin{equation*}
\ka:X\setm Z\ra M
\end{equation*}

If $\Lb\in\Pic^{-e}(C)$ denote by $\ka_{\Lb}$ the restriction of $\ka$ to
the fiber of $p:X\ra\Pic^{-e}(C)$ at $\{\Lb\}$
\begin{equation*}
\ka_{\Lb}:\PP(V_{\Lb})\setm Z_{\Lb}\ra M
\end{equation*}

Note that $\dim\PP(V_{\Lb})=2e+g-1$, $\dim X=2e+2g-1$, $\dim M=3g-3$.

\begin{lemma}\label{ka_injective}
If $e\leq(\frac{g}{2}-1)$, for $\Lb\in\Pic^{-e}(C)$ general, the
morphism $\ka_{\Lb}$ is injective.
\end{lemma}

\bp
Assume that there are two extensions in
$V_{\Lb}=\Ext_C^1(\Lb^{-1}\otimes\xi,\Lb)$ with the same middle term.
If the two extensions correspond to different elements in $\PP(V_{\Lb})$,
then there is a non-zero morphism $\Lb\ra\Lb^{-1}\otimes\xi$,
i.e., $\H^0(\Lb^{-2}\otimes\xi)\neq0$. But if $e\leq(\frac{g}{2}-1)$ and
$\Lb\in\Pic^{-e}(C)$ is general, then $\H^0(\Lb^{-2}\otimes\xi)=0$ (see
\cite{ACGH}, IV.4.5) and our conclusion follows.
\ep

\begin{cor}
If $\Lb\in\Pic^0(C)$ is general, the morphism $\ka_{\Lb}$
gives a linear embedding of $\PP(V_{\Lb})\cong\PP^{g-1}$ into $M$.
\end{cor}

It follows from \cite{B} (Prop. 4.5) that if $e\geq(g-1)$, for any
$\Lb\in\Pic^{-e}(C)$, the morphism $\ka_{\Lb}$ is dominant.
Note that the bound for $e$ is optimal.



\section{The nice component of the space of rational curves of odd degree}

If $k$ is a positive integer, let $\Hom_k(\PP^1,M)$ be the Hilbert
scheme of morphisms $\PP^1\ra M$ of degree $k$ \cite{K}. The scheme
$\Hom_k(\PP^1,M)$ is an open in a Hilbert scheme of rational curves in
$\PP^1\times M$ (the graph of $f$ is used to embed $\PP^1$ into
$\PP^1\times M$).


\subsection{The degree of the line bundle $\ka_{\Lb}^*\Th$}
For $\Lb\in\Pic^{-e}(C)$ consider the morphism
$\ka_{\Lb}:\PP(V_{\Lb})\setm Z_{\Lb}\ra M$.
By Lemma \ref{unstable_locus_1_local}, the codimension of $Z_{\Lb}$ in
$\PP(V_{\Lb})$ is at least $g\geq2$.
We compute the degree of the line bundle $\ka_{\Lb}^*\Th$, using the
following more general result.

\begin{lemma}\label{degree}
Let $S$ be a smooth variety and $\F$ a rank $2$ vector
bundle on $S\times C$, such that, for any $s\in S$, the restriction
$\F_s$ of $\F$ to $\{s\}\times C$ is stable, of determinant $\xi$, so that
$\F$ induces a morphism $f:S\ra M$. Denote by $\pi_1$ and $\pi_2$ the two
projections from $S\times C$. We have:
\begin{equation*}
c_1(f^*(\Th))={\pi_1}_*(2c_2(\F)-\frac{c_1(\F)^2}{2})
\end{equation*}
\end{lemma}

\bp
Denote by $\ad(\F)$ the subbundle of $\EEnd(\F)\cong\F\otimes\F^*$ of
endomorphisms of trace zero. There is a split exact sequence:
\begin{equation*}
0\ra\ad(\F)\ra\F\otimes\F^*\ra O\ra0
\end{equation*}

Let $T_M$ be the tangent bundle of $M$.
We have from \cite{N2} that:
\begin{equation}\label{f^*T}
f^*(T_M)\cong\R^1{\pi_1}_*(\ad(\F))
\end{equation}

If $\HH$ is a bundle on $S\times C$, by the Grothendieck-Riemann-Roch
formula, we have the following relations in $A(S)\otimes\QQ$:
\begin{gather*}
ch({\pi_1}_!(\HH))={\pi_1}_*(ch(\HH).td({\pi_2}^*T_C))\\
c_1({\pi_1}_!(\HH))=\frac{1}{2}
{\pi_1}_*(c_1^2(\HH)-2c_2(\HH)-c_1(\HH).{\pi_2}^*c_1(K_C))
\end{gather*}
We apply this to $\HH=\ad(\F)$. We have ${\pi_1}_*(\ad(\F))=0$,
$c_1(\ad(\F))=0$ and $$c_2(\ad(\F))=c_2(\F\otimes\F^*)=4c_2(\F)-c_1(\F)^2.$$
Combining all the relations, one gets the formula:
\begin{equation}\label{GRR_relation}
c_1(\R^1{\pi_1}_*(\ad(\F)))={\pi_1}_*(4c_2(\F)-c_1(\F)^2)
\end{equation}
As  $K_M\cong\Th^{-2}$, the result follows now from
(\ref{f^*T}) and (\ref{GRR_relation}).
\ep

\begin{cor}\label{deg_1}
We have $\ka_{\Lb}^*(\Th)\cong O(2e+1)$.
\end{cor}

\bp
Consider the restriction of the universal sequence
(\ref{univ_ext_1}) to $\PP(V_{\Lb})\times C$:
\begin{equation*}
0\ra\nu_2^*\Lb(1)\ra\G_{\Lb}\ra\nu_2^*(\Lb^{-1}\otimes\xi)\ra0.
\end{equation*}
where $\G_{\Lb}$ is the restriction of the bundle $\G$ to
$\PP(V_{\Lb})\times C$. If we let $\{H\}\in\AA(\PP(V_{\Lb}))$ be the class
of a hyperplane in $\PP(V_{\Lb})$, we have:
\begin{equation*}
c_1(\G_{\Lb})=\{H\}\times C + \PP(V_{\Lb})\times c_1(\xi),\quad
c_2(\G_{\Lb})=\{H\}\times (c_1(\xi)-c_1(\Lb))
\end{equation*}
The result follows from Lemma \ref{degree}.
\ep


\subsection{The nice component $\M$ of $\Hom_k(\PP^1, M)$ for $k=2e+1$}

\begin{lemma}\label{unobstr_1}
Let $f:\PP^1\ra M$ be the composition
\begin{equation}\label{comp_1}
\begin{CD}
\PP^1 @>h>>\PP(V_{\Lb})\setm Z_{\Lb} @>\ka_{\Lb}>> M
\end{CD}
\end{equation}
where $h$ is such that $h^*O(1)\cong O(1)$. Then $f$ is
an unobstructed point of $\Hom_k(\PP^1, M)$.
\end{lemma}

\bp
Note that by Corollary \ref{deg_1}, we have that $f^*\Th=O(k)$, so $f$
is a point of $\Hom_k(\PP^1, M)$.
Consider the exact sequence obtained by pulling back the
universal extension (\ref{univ_ext_1}) to $\PP^1\times C$ by $(g\times\id_C)$:
\begin{equation}\label{seq_1}
0\ra O(1)\bo\Lb\ra\F\ra O\bo(\Lb^{-1}\otimes\xi)\ra0.
\end{equation}
The bundle $\F$ induces the morphism $f:\PP^1\ra M$.
It follows that if $\eta$ is the generic point of $C$, then
$\F_{\eta}\cong O(1)\oplus O$. From the more general Lemma
\ref{useful_result}, it follows that $f$ is an unobstructed
point in $\Hom_k(\PP^1,M)$.
\ep

\begin{lemma}\label{useful_result}
Let $f:\PP^1\ra M$ be a morphism given by the vector bundle $\F$
on $\PP^1\times C$ with the property that, if $\eta$ is the generic point of
$C$, then the bundle $\F_{\eta}$ on $\PP^1_{\eta}$ is balanced, i.e.,
it splits as $O(a)\oplus O(b)$, for some integers $a$ and $b$ with
$|b-a|\geq1$. Then $\H^1(\PP^1,f^*T_{M})=0$.
\end{lemma}

\bp
We have from \cite{N2} that $f^*T_{M}\cong\R^1{p_1}_*(\ad(\F))$ and
$\R^i{p_1}_*(\ad(\F))=0$, for any $i\neq1$ (where $p_1$ and $p_2$ are the
two projections from $\PP^1\times C$). There are two Leray spectral
sequences:
\begin{gather*}
\H^i(\PP^1,\R^j{p_1}_*(\ad(\F)))\Lra
\H^{i+j}(\PP^1\times C,\ad(\F))\\
\H^i(C,\R^j{p_2}_*(\ad(\F)))\Lra
\H^{i+j}(\PP^1\times C,\ad(\F))
\end{gather*}
Then, as for $i\geq2$ we have $\R^i{p_2}_*(\ad(\F))=0$, it follows that
\begin{equation}\label{h^1}
\H^1(\PP^1,\R^1{p_1}_*(\ad(\F)))\cong
\H^2(\PP^1\times C,\ad(\F))\cong\H^1(\R^1{p_2}_*(\ad(\F)))
\end{equation}
The sheaf  $\R^1{p_2}_*(\ad(\F))$ is supported at a finite number
of points of $C$, as we have:
\begin{equation*}
\R^1{p_2}_*(\ad(\F))_{\eta}\cong\H^1(\PP^1_{\eta},\ad(\F_{\eta}))=0
\end{equation*}
It follows that $\H^1(\R^1{p_2}_*(\ad(\F))=0$ and the result follows from
(\ref{h^1}).
\ep

\begin{thm}\label{nice_comp_odd}
If $k=2e+1$ is an odd positive integer, there is an irreducible component
$\M$ of $\Hom_k(\PP^1,M)$, having the expected dimension $2k+3g-3$, with
the following properties:
\bi
\item[i)] A general point $f\in\M$ is unobstructed and it is
obtained as a composition:
\begin{equation*}
\begin{CD}
\PP^1@>h>>\PP(V_{\Lb})\setm Z_{\Lb}@>{\ka_{\Lb}}>>M
\end{CD}
\end{equation*}
where $\Lb\in\Pic^{-e}(C)$ and $h$ is a morphism such that $h^*O(1)\cong O(1)$
\item[ii)] The MRC fibration of $\M$ is given by a dominant rational map:
\begin{equation*}
\M\dra\Pic^{-e}(C)
\end{equation*}
which sends a point $f\in\M$ as in $i)$ to $\Lb\in\Pic^{-e}(C)$.
\ei
\end{thm}

\bp
Abusing notations, we let $\Hom_1(\PP^1,X\setm Z)$ to be the Hilbert
scheme of morphisms $\PP^1\ra X\setm Z$ that are contained
in the fibers of the projective bundle $p:X\ra\Pic^{-e}(C)$ and have
degree $1$ with respect to the relatively ample line bundle $O_X(1)$
(see \cite{K}). There is a canonical morphism:
\begin{equation}\label{pi}
\pi:\Hom_1(\PP^1,X\setm Z)\ra\Pic^{-e}(C)
\end{equation}
whose fiber at $\Lb\in\Pic^{-e}(C)$ is
$\Hom_1(\PP^1,\PP(V_{\Lb})\setm Z_{\Lb})$.
(Note that, as $\Pic^{-e}(C)$ is an abelian variety, any morphism
$\PP^1\ra X$ must map to a point in $\Pic^{-e}(C)$.)

One has to notice that a morphism $\PP^1\ra\PP^r$ is unobstructed.
We have $\Hom_1(\PP^1,\PP^r)\cong U$, where $U\subset\PP^{2r+1}$ is a dense
open set. It follows that $\Hom_1(\PP^1,X\setm Z)$ is a smooth,
irreducible variety of the expected dimension:
\begin{equation*}
\dim\Hom_1(\PP^1,X)=g+2(2e+g-1)+1=2k+3g-3
\end{equation*}
Then $\ka:X\setm Z\ra M$ induces a morphism between the corresponding
schemes of morphisms:
\begin{equation}\label{psi}
\psi:\Hom_1(\PP^1,X\setm Z)\ra\Hom_k(\PP^1,M)
\end{equation}

We prove that the morphism $\psi$ is injective. Let $h$ and $h'$ be morphisms
$\PP^1\ra X\setm Z$, such that $\ka\circ h=\ka\circ h'$. Assume that $h$,
respectively $h'$, have image in the fiber of $p$ over
$\Lb\in\Pic^{-e}(C)$, respectively $\Lb'\in\Pic^{-e}(C)$. We have that the
following two compositions are equal:
\begin{equation*}
\begin{CD}
\PP^1 @>h>>\PP(V_{\Lb})\setm Z_{\Lb} @>\ka_{\Lb}>> M\quad\hbox{and}
\quad\PP^1 @>h'>>\PP(V_{\Lb'})\setm Z_{\Lb'} @> \ka_{\Lb'}>> M
\end{CD}
\end{equation*}
As in (\ref{seq_1}), there are two exact sequences on $\PP^1\times C$:
\begin{equation}\label{seq_1_F}
0\ra O(1)\bo\Lb\ra\F\ra O\bo(\Lb^{-1}\otimes\xi)\ra0.
\end{equation}
\begin{equation}\label{seq_1_F'}
0\ra O(1)\bo\Lb'\ra\F'\ra O\bo({\Lb'}^{-1}\otimes\xi)\ra0.
\end{equation}

Since $\F$ and $\F'$ induce the same morphism $\PP^1\ra M$,
there is an integer $m$ such that $\F'\cong\F\otimes p_1^*O(m)$
(without loss of generality, we may assume $m\geq0$). As there are no
non-zero morphisms $O(1)\ra O(-m)$, there is a commutative diagram:
\begin{equation}\label{diagram_1_F_F'}
\begin{CD}
0@>>>O(1)\bo\Lb@>>>\F@>>>O\bo(\Lb^{-1}\otimes\xi)@>>> 0\\
@VVV @VVV  @VV{\cong}V @VVV @VVV\\
0@>>>O(1-m)\bo\Lb'@>>>\F'(-m)@>>>O(-m)\bo({\Lb'}^{-1}\otimes\xi)@>>>0
\end{CD}
\end{equation}

Since a surjective  morphism of line bundles is an isomorphism, it follows
that all
the vertical arrows in (\ref{diagram_1_F_F'}) are isomorphisms. Hence, $m=0$,
$\Lb\cong\Lb'$ and $\F\cong\F'$. As endomorphisms of line bundles are given
by multiplication with scalars, the extensions
(\ref{seq_1_F}) and (\ref{seq_1_F'}) are scalar multiples of each other.
Then we must have $h=h'$, as $h$ sends a  point $p\in\PP^1$ to the class
in $\PP(V_{\Lb})$ of the extension obtained by restricting (\ref{seq_1_F}) to
$\{p\}\times C$
\begin{equation*}
0\ra\Lb\ra\F_p\ra \Lb^{-1}\otimes\xi\ra0
\end{equation*}

It follows that $\psi$ is injective. Let $\M$ be the closure of the image of
the morphism $\psi$ with the reduced structure.
We have that $\M$ is an irreducible closed subvariety of $\Hom_k(\PP^1,M)$,
birational to $\Hom_1(\PP^1, X\setm Z)$, which has the expected dimension
$2k+3g-3$. By Lemma \ref{unobstr_1}, a general point $f\in\M$ is unobstructed.
It follows that $\M$ is the unique irreducible component of $\Hom_k(\PP^1,M)$
containing the point $f$. This proves i).
\ps

For ii), consider the morphism $\pi:\Hom_1(\PP^1, X\setm Z)\ra\Pic^{-e}(C)$
of (\ref{pi}). Since $\M$ is birational to $\Hom_1(\PP^1, X\setm Z)$ ,
there is a rational map:
\begin{equation*}
\rho:\M\dra\Pic^{-e}(C)
\end{equation*}
We claim that this is the MRC fibration of $\M$. As $\Pic^{-e}(C)$ is an
abelian variety, so it does not contain any rational curves, it is enough to
prove that the fibers of $\rho$ are rationally connected (or, equivalently,
the fibers of $\pi$ are rationally connected). But the fiber of $\pi$ at a
point $\Lb$ is $\Hom_1(\PP^1, \PP(V_{\Lb})\setm Z_{\Lb})$, which is a rational
variety, and this proves ii).
\ep


\subsection{The subvarieties $\S(e,n)\subset\Hom_k(\PP^1,M)$}\label{S(e,n)}
If $n\geq1$ is any integer, we let $k=n(2e+1)$. Using the same arguments as
in the proof of Theorem \ref{nice_comp_odd}, one can still prove that there
is an injective morphism:
\begin{equation*}
\psi:\Hom_n(\PP^1,X\setm Z)\ra\Hom_k(\PP^1,M)
\end{equation*}

It follows that there are irreducible subvarieties
$\S(e,n)\subset\Hom_k(\PP^1,M)$ such that:
\bi
\item[i) ]A general point $f\in\S(e,n)$ is obtained as a composition
\begin{equation}\label{comp_1'}
\begin{CD}
\PP^1@>h>>\PP(V_{\Lb})\setm Z_{\Lb}@>{\ka_{\Lb}}>>M_{\xi}
\end{CD}
\end{equation}
where $\Lb\in\Pic^{-e}(C)$ and $h$ is a morphism such that
$h^*O(1)\cong O(n)$
\item[ii)]The MRC fibration of $\S(e,n)$ is given by a rational map:
\begin{equation*}
\S(e,n)\dra\Pic^{-e}(C)
\end{equation*}
which sends a general point $f\in\M$ as in $i)$ to $\Lb\in\Pic^{-e}(C)$.
\ei
Note that for $n=1$ the subvariety $\S(k,1)$ is the nice component $\M$.


\subsection*{Remark}
If $n>1$, we cannot apply the Lemma \ref{useful_result}
to prove that a general point $f\in\S(e,n)$ is unobstructed.
Note that the dimension of $\S(e,n)$ is equal to:
\begin{equation*}
\dim\Hom_n(\PP^1,X\setm Z)=g+(2e+g)(n+1)
\end{equation*}

For $n>1$ we have that $\dim\S(e,n)\geq(2k+3g-3)$ if and only if
$e\leq\frac{g}{2}-1$. As any irreducible component has at
least the expected dimension, it follows that if $n>1$ and
$e>\frac{g}{2}-1$, then $\S(e,n)$ is not an irreducible component.

If $e<\frac{g}{2}-1$, then $\dim\S(e,n)>(2k+3g-3)$. It follows that
any point of $\S(e,n)$ is obstructed.
With our methods so far, we are not able to prove
that $\S(e,n)$ is an irreducible component. We will prove in \cite{C} that
this is indeed the case.

If $e=\frac{g}{2}-1$ ($g$ even), note that $k$ must be divisible by $(g-1)$.
We prove in \ref{almost_nice} that this is another irreducible component.


\section{Free rational curves on $M$  and the almost nice component}

\subsection{Free rational curves}

We recall some terminology. If $\G\cong\oplus_i O(a_i)$ is a vector bundle
on $\PP^1$, we say that $\G$ is \emph{positive} (resp.
\emph{non-negative}), if for all $i$, we have $a_i>0$ (resp. $a_i\geq0$).
We say that a rational curve $f:\PP^1\ra M$  is \emph{very free}
(resp. \emph{free}), if the bundle $f^*T_M$ is positive (resp. non-negative).

\begin{lemma}\label{free_1}
Let $f:\PP^1\ra M$ be a composition as in (\ref{comp_1'}).
For $\Lb\in\Pic^{-e}(C)$ general, we distinguish the following cases:
\bi
\item[i. ]If $e<\frac{g}{2}-1$, then
$f^*T_M\cong O(-n)^{g-2e-2}\oplus O^g\oplus\A$, for some $\A$ positive
\item[ii. ]If $e=\frac{g}{2}-1$ ($g$ even), then
$f^*T_M\cong O^g\oplus\A$, for some $\A$ positive
\item[iii. ]If $e\geq\frac{g}{2}-1$, then $f$ is a free curve, with at most
$g$ trivial components
\ei
Moreover, if $e\geq2g-1$, for any $\Lb\in\Pic^{-e}(C)$, we have that $f$
is a very free curve.
\end{lemma}

\bp
Recall that $f:\PP^1\ra M$ is induced by the vector bundle $\F$ on
$\PP^1\times C$,
then $f^*T_M=\R^1{p_1}_*\ad(\F)$. As
$\F\otimes\F^*\cong\ad(\F)\oplus O$,
it follows that $\R^1{p_1}_*(\F\otimes\F^*)\cong f^*T_M\oplus O^g$.
We denote $$\T=\R^1{p_1}_*(\F\otimes\F^*).$$

We compute $\T$ for $f$ as
in (\ref{comp_1'}). There is an exact sequence on $\PP^1\times C$:
\begin{equation}\label{seq_F_1}
0\ra O(n)\bo\Lb\ra\F\ra O\bo(\Lb^{-1}\otimes\xi)\ra0
\end{equation}

Denote
\begin{equation*}
\S'=(O\bo(\Lb\otimes\xi^{-1}))\otimes\F,\quad\quad
\S''=(O(-n)\bo\Lb^{-1})\otimes\F
\end{equation*}

\ps

By tensoring (\ref{seq_F_1}) with
$O\bo(\Lb\otimes\xi^{-1})$, resp. with $O(-n)\bo\Lb^{-1}$, we get
exact sequences:
\begin{equation}\label{II}
0\ra O(n)\bo(\Lb^2\otimes\xi^{-1})\ra\S'\ra O\ra0
\end{equation}
\begin{equation}\label{III}
0\ra O\ra\S''\ra O(-n)\bo(\Lb^{-2}\otimes\xi)\ra0
\end{equation}

By dualizing  (\ref{seq_F_1}) and tensoring with $\F$, we get an exact 
sequence:
\begin{equation}\label{I}
0\ra\S'\otimes\F\ra\F\otimes\F^*\ra\S''\ra0
\end{equation}

Note that $\H^0(C,\S'_p)\cong\Hom(\Lb^{-1}\otimes\xi,\F_p)$.
Since $\F_p$ is stable, of degree $1$, for any $p\in\PP^1$, there are no
non-zero morphisms $\Lb^{-1}\otimes\xi\ra\F_p$.
It follows that ${p_1}_*(\S')=0$.

\ps

Applying ${p_1}_*$ to (\ref{II}) and using that ${p_1}_*(\S')=0$, we have
an exact sequence:
\begin{equation}\label{I'}
0\ra O\ra O(n)\otimes\H^1(C,\Lb^2\otimes\xi^{-1})\ra\R^1{p_1}_*(\S')
\ra O^g\ra0
\end{equation}

Applying ${p_1}_*$ to (\ref{III}) we have another sequence
\begin{equation}\label{II'}
\begin{CD}
0\ra O\ra{p_1}_*(\S'')\ra O(-n)\otimes\H^0(C,\Lb^{-2}\otimes\xi)\ra O^g\ra\\
\ra\R^1{p_1}_*(\S'')\ra O(-n)\otimes\H^1(C,\Lb^{-2}\otimes\xi)\ra0
\end{CD}
\end{equation}

Applying ${p_1}_*$ to (\ref{I}), and using that
${p_1}_*(\F\otimes\F^*)\cong O$, we have a sequence
\begin{equation}\label{T}
0\ra O\ra{p_1}_*(\S'')\ra\R^1{p_1}_*(\S')\ra\T\ra\R^1{p_1}_*(\S'')\ra0
\end{equation}

From (\ref{I'}), it follows that there is a positive vector bundle $\A$
such that
\begin{equation}\label{A}
\R^1{p_1}_*(\S')\cong O^g\oplus\A
\end{equation}

\ps

If $e\geq\frac{g}{2}-1$ and $\Lb\in\Pic^{-e}(C)$ is general, then
$\H^1(C,\Lb^{-2}\otimes\xi)=0$ (see \cite{ACGH}, IV.4.5).
In this case, we have from (\ref{II'}) and (\ref{T}) that
$\R^1{p_1}_*(\S'')$ and $\T$ are non-negative. It follows that $f^*T_M$ is
non-negative. Note that from (\ref{II'}) we have that $\R^1{p_1}_*(\S'')$
has at most $g$ trivial components. It follows from
(\ref{A}) and (\ref{T}) that $\T$ has at most $2g$ trivial components.
Hence, $f^*T_M$ has at most $g$ trivial components.
This proves iii).

\ps

If $e=\frac{g}{2}-1$, then for $\Lb\in\Pic^{-e}(C)$ general,
we have  $\H^0(C,\Lb^{-2}\otimes\xi)=0$.
It follows from (\ref{II'}) that ${p_1}_*(\S'')\cong O$ and
$\R^1{p_1}_*(\S'')\cong O^g$. From (\ref{T}), there is an exact sequence:
\begin{equation*}
0\ra\R^1{p_1}_*(\S')\ra\T\ra\R^1{p_1}_*(\S'')\ra0
\end{equation*}

As $\T\cong f^*T_M\oplus O^g$, it follows using (\ref{A})
that $f^*T_M\cong O^g\oplus\A$. This proves ii).

\ps

We prove now i). If $e<\frac{g}{2}-1$ and $\Lb\in\Pic^{-e}(C)$ is general,
then $\H^0(C,\Lb^{-2}\otimes\xi)=0$. By Riemann-Roch,
$\h^1(C,\Lb^{-2}\otimes\xi)=g-2e-2$. It follows from  (\ref{II'}) that
${p_1}_*(\S'')\cong O$ and there is an exact sequence:
\begin{equation*}
0\ra O^g\ra\R^1{p_1}_*(\S'')\ra O(-n)^{g-2e-2}\ra0
\end{equation*}

Then $\R^1{p_1}_*(\S'')\cong O^g\oplus O(-n)^{g-2e-2}$.
We have from (\ref{T}) that there is an exact sequence:
\begin{equation}\label{T_2}
0\ra\R^1{p_1}_*(\S')\ra\T\ra\R^1{p_1}_*(\S'')\ra0
\end{equation}

As $\R^1{p_1}_*(\S'')\cong O^g\oplus O(-n)^{g-2e-2}$,
it follows from (\ref{A}) that the sequence (\ref{T_2}) is split.
Since $\T\cong f^*T_M\oplus O^g$, we have:
\begin{equation*}
f^*T_M\cong O(-n)^{g-2e-2}\oplus O^g\oplus\A
\end{equation*}
This proves i).

\ps

Assume $e\geq2g-1$.
We have $\deg(K^{-1}_C\otimes\Lb^{-1})\geq0$.
Since for any $p\in\PP^1$, the bundle $\F^*_p$ stable, of degree $-1$,
there are no non-zero morphisms
$K_C^{-1}\otimes\Lb^{-1}\ra\F^*_p$. It follows that
\begin{equation*}
\H^0(C,K_C\otimes\Lb\otimes\F^*_p)\cong\Hom(K^{-1}_C\otimes\Lb^{-1},\F_p)=0
\end{equation*}

By Serre duality, we have $\H^1(C,\Lb^{-1}\otimes\F_p)=0$ for any
$p\in\PP^1$. Then $\R^1{p_1}_*(\S'')=0$. From (\ref{T}) we have
an exact sequence:
\begin{equation*}
0\ra O\ra{p_1}_*(\S'')\ra\R^1{p_1}_*(\S')\ra\T\ra0
\end{equation*}
As $\R^1{p_1}_*(\S')=O^g\oplus\A$ ($\A$ positive) and
$\T=O^g\oplus f^*T_M$, it follows that $f^*T_M$ is positive.
\ep

\begin{cor}\label{very_free_odd}
Let $k$ be an odd integer and let $\M\subset\Hom_k(\PP^1,M)$ be the nice
component. We distinguish the following cases:
\bi
\item[i) ]If $k\leq(g-1)$, then for general $f\in\M$, we have
$$f^*T_M\cong O(-1)^{g-k-1}\oplus O^g\oplus O(1)^{k+g-3}\oplus O(2)$$
\item[ii) ]If $k\geq(g-1)$, then a general $f\in\M$ is a free curve
with at most $g$ trivial components
\item[iii) ]If $k\geq(2g-1)$, then a general $f\in\M$ is a very free curve
\ei
\end{cor}

\bp
Parts i) and ii) are immediate consequences of Lemma \ref{free_1} for the
case when $n=1$ and $k=2e+1$. Note that in this case, the bundle $\A$
from (\ref{A}) can be computed from (\ref{I'}) to be
$O(1)^{2e+g-2}\oplus O(2)$.

For part iii), note that we have from Lemma \ref{free_1} that if
$k\geq(4g-1)$, then a general $f\in\M$ is a very free curve.
We may improve this bound by observing that if
$k\geq(2g-1)$, then the morphism $\ka_{\Lb}$
is dominant (\cite{B}, Prop.4.5).
For general $\Lb\in\Pic^{-e}(C)$, we have by Lemma
\ref{free_1} that for any $h:\PP^1\ra\PP(V_{\Lb})\setm Z_{\Lb}$,
if we let $f=\ka_{\Lb}\circ h$, then $f$ is a free curve.
Since any degree $1$ morphism $h$ is very free (as
$h^*T_{\PP(V_{\Lb})}\cong O(1)^{2e+g-2}\oplus O(2)$), it follows
from Lemma \ref{very_free} that the morphism $f$ is very free.
\ep

\begin{lemma}\label{very_free}
Let $X$ and $Y$ be smooth varieties and let $\phi:X\ra Y$ be a
dominant morphism. Let $h:\PP^1\ra X$ be a rational curve
intersecting the smooth locus of the morphism $\phi$. Let
$f=\phi\circ h$ and assume that $f$ is a free curve and that $h$
is very free. Then $f$ must be very free.
\end{lemma}

\bp Denote with $u$ the morphism $h^*T_X\ra f^*T_Y$ of vector
bundles on $\PP^1$. Assume that $f$ is not very free, hence, the
bundle $f^*T_Y$ is non-negative with $r>0$ trivial components. As
$h^*T_X$ is positive, it follows that $O^r$ is a subsheaf of the
cokernel of $u$. But since $h$ intersects the smooth locus of
$\phi$, the morphism $u$ is generically surjective, so its
cokernel is supported at a finite number of points and we get a
contradiction. \ep


\subsection{The almost nice component of $\Hom_k(\PP^1,M)$}
\label{almost_nice}

\begin{thm}
If $g$ is even and $k>0$ is an integer that is divisible by
$(g-1)$, there is an irreducible component
$\M'\subset\Hom_k(\PP^1,M)$, having the expected dimension
$(2k+3g-3)$ and with the following properties: \bi \item[i) ]A
general point $f\in\M'$ is obtained as a composition:
\begin{equation*}
\begin{CD}
\PP^1@>h>>\PP(V_{\Lb})\setm Z_{\Lb}@>{\ka_{\Lb}}>>M
\end{CD}
\end{equation*}
where $\Lb\in\Pic^{1-\frac{g}{2}}(C)$ and $h$ is a morphism such that
$h^*O(1)\cong O(\frac{k}{g-1})$
\item[ii)] The MRC fibration of $\M'$ is given by a rational map:
\begin{equation*}
\M\dra\Pic^{1-\frac{g}{2}}(C)
\end{equation*}
which sends a point $f\in\M$ as in $i)$ to
$\Lb\in\Pic^{1-\frac{g}{2}}(C)$. \item[iii) ]A general point
$f\in\M'$ has the property that $f^*T_M\cong O^g\oplus\A$, for
some positive vector bundle $\A$ \ei
\end{thm}

\bp Let $\M'$ be the irreducible subvariety
$\S(e,n)\subset\Hom_k(\PP^1,M)$ constructed in \ref{S(e,n)}, for
$e=\frac{g}{2}-1$ and $n=\frac{k}{g-1}$. We already have i) and
ii) from \ref{S(e,n)} and iii) from Lemma \ref{free_1}. We only
have to prove that $\M'$ is an irreducible component.

We noted in \ref{S(e,n)} that the dimension of $\M'$ is equal to the
expected dimension. By Lemma \ref{free_1}, a general point $f\in\M'$ is
unobstructed. It follows that $\M'$ is the unique irreducible component of
$\Hom_k(\PP^1,M)$ containing the point $f$.
\ep

\begin{lemma}
If $g\geq3$ the general element in the almost nice component is an
embbeded curve.
If $g=2$, the general element in the almost nice component is a morphism
$\PP^1\ra M$ that is $k$-to-$1$ onto  a line in $M$.
\end{lemma}

\bp If $e=\frac{g}{2}-1$, and $\Lb$ is general in $\Pic^{-e}(C)$,
the morphism $\ka_{\Lb}:\PP(V_{\Lb})\setm Z_{\Lb}\ra M$ is
injective, by Lemma \ref{ka_injective}. Since
$\dim\PP(V_{\Lb})=2g-3$, a general morphism $\PP^1\ra\PP(V_{\Lb})$
of degree $\frac{k}{g-1}$, will give by composition with
$\ka_{\Lb}$ an injective morphism $\PP^1\ra M$ if and only if
$g\geq3$.  If $g=2$ then $e=0$ and by Lemma \ref{ka_injective} it
follows that $\ka_{\Lb}$ gives a linear embedding of $\PP^1$ into
$M$, and the result follows.

\ep


\section{Extensions of skyscraper sheaves by rank $2$ vector bundles}

Let $e\geq1$ be an integer, $D\in\Sym^e(C)$ and let $\E$ be a rank $2$,
semistable vector bundle on $C$, with $\det(\E)\cong\xi(-D)$.
Consider extensions:
\begin{equation}\tag{$**$}
0\ra\E\ra\E'\ra O_D\ra0
\end{equation}
Such extensions are classified by the vector space:
\begin{equation}\label{V_D,E}
V_{D,\E}:=\Ext^1_C(O_D, \E)\cong\H^0(C, \E(D)_{|D})
\end{equation}

Then the vector space $V_{D,\E}$ has dimension $2e$. The
isomorphism classes of non-trivial extensions as above are
parametrized by the projective space $\PP(V_{D,\E})$.

\ps

Let $M_{1-e}$ be the coarse moduli scheme of rank $2$, semistable vector
bundles on $C$, with determinant of degree $(1-e)$. Let
the determinant map be $\det:M_{1-e}\ra\Pic^{1-e}(C)$.
Let $M^s_{1-e}\subset M_{1-e}$ be the locus of stable vector bundles.
We now define a parameter space for extensions $(**)$, while we let
$D\in\Sym^e(C)$ and $\E\in M^s_{1-e}$ vary, with the condition that
$\det(\E)\cong\xi(-D)$.

\ps

The moduli space $M_{1-e}$ is the geometric quotient of a smooth,
quasiprojective variety $\MM_{1-e}$, by the action of an algebraic group
$\PGL(r)$, for some integer $r$.
Let $\tau:\MM_{1-e}\ra M_{1-e}$ be the quotient map and
let $\MM^s_{1-e}=\tau^{-1}(\MM^s_{1-e})$.
Note that there is a {Poincar\'e} bundle $\V$ on $\MM_{1-e}\times C$.

\ps

Define $\th:\Pic^{1-e}(C)\ra\Pic^e(C)$ to be the morphism that sends $\Lb$ to
$\xi\otimes\Lb^{-1}$. Consider the composition
$\th\circ(\det)\circ\tau:\MM^s_{1-e}\ra\Pic^e(C)$.
We define  $S'$ to be the fiber product
\begin{equation}\label{S'}
S'=\MM^s_{1-e}\times_{\Pic^e(C)}\Sym^e(C)
\end{equation}

Note that $S'$ is a smooth variety of dimension $(3g-3)+e+\dim\PGL(r)$.

\ps

Let $\V_{S'}$ be the pull-back of the bundle $\V$ to $S'\times C$
and let $\D_{S'}\subset S'\times C$ be the universal scheme coming
from $D\in\Sym^e(C)$. Define on $S'$ the relative extension sheaf
\begin{equation*}
\S=\EExt^1_{S'\times C|S'}(O_{\D_{S'}},\V_{S'})
\end{equation*}

Note that $\S$ is a locally free sheaf. Let  $W'$ be the
projective bundle $\PP(\S)$ and let
\begin{equation}\label{W'}
p:W'\ra S'
\end{equation}
be the corresponding map. If $\mu_1,\mu_2$ are the two projections
from $W'\times C$, there is a universal extension on $W'\times C$:
\begin{equation}\label{univ_ext_2}
0\ra \mu_1^*O_{W'}\otimes p^*\V_{S'}\ra\J\ra p^*O_{\D_{S'}}\ra0.
\end{equation}


\subsection{The locus of not locally free extensions}

Consider the space of extensions $W'$. Let $\Ga\subset W'$ be the locus of
those $w\in W'$, for which the bundle $\J_w$ of (\ref{univ_ext_2}) is not
locally free. The following lemma gives an estimate of its codimension.

\begin{lemma}\label{not_lolally_free}
For $D\in\Sym^e(C)$ and $\E$ a rank $2$ vector bundle such that
$\det(\E)\cong\xi(-D)$, there is a closed subvariety
$\Ga_{D,\E}\subset\PP(V_{D,\E})$, corresponding to
not locally free extensions in $\PP(V_{D,\E})$. If $e\geq2$, then
$\Ga_{D,\E}$ consists of a union of $e$ linear subspaces of codimension $2$.
If $e=1$, then $\Ga_{D,\E}=\eset$.
\end{lemma}

\bp Denote $V=V_{D,\E}$. Let $D$ consist of points
$y_1,\ldots,y_e$ on $C$ (not necessarily distinct). For
$i\in\{1,\ldots,e\}$, consider the $2$-dimensional vector spaces
$V_i=\Ext^1_C(O_{y_i},\E)$. There is an isomorphism
$V\cong\oplus_{i=1}^eV_i$.
If $(**)$ corresponds to $v=(v_1,\ldots,v_e)\in V$, then it
follows that $\E'$ is not locally free if and only if $v_i=0$ for
some $i$. Let $\Ga=\Ga_{D,\E}$. If $e=1$, then $\Ga=\eset$. If
$e\geq2$, then $\Ga\subset\PP(V)$ is the union of the linear
spaces $v_i=0$. \ep

\begin{cor}
There is a closed subvariety $\Ga\subset W'$, corresponding to extensions
in $W'$ that are not locally free, and all the irreducible components of
$\Ga$ have codimension at least $2$ if $e\geq2$. If $e=1$, then $\Ga=\eset$.
\end{cor}

\subsection{Action of $\Aut(O_D)$ on $\PP(V_{D,\E})$}\label{G_action}

Consider an extension
\begin{equation*}
0\ra\E\ra\E'\ra O_D\ra0
\end{equation*}

Clearly, the group $\Aut(O_D)$ acts on $V_{D,\E}$, as well as on
$\PP(V_{D,\E})$.

\ps

Let $D$ consist of distinct points $p_1,\ldots,p_e$. If we let
$V=V_{D,\E}$ and $x_i, y_i$ be coordinates on
$V_i=\Ext(O_{p_i},\E)$, then the group $\Aut(O_D)\cong\GGG_m^e$
acts on $\PP(V)$ by:
\begin{equation*}
(\la_1,\ldots,\la_e).(x_1,y_1,\ldots,x_e,y_e)\mapsto
({\la_1}x_1,{\la_1}y_1,\ldots,{\la_e}x_e,{\la_e}y_e)
\end{equation*}


\begin{lemma}\label{fibers}
Let $D$ consist of distinct points $p_1,\ldots,p_e$. The orbits
for the action of $\Aut(O_D)$ on $\PP(V_{D,\E})$ are
$(e-1)$-dimensional linear subspaces of $\PP(V_{D,\E})$.
\end{lemma}

\bp
Consider an extension corresponding to a point $p\in\PP(V)$:
\begin{equation*}
0\ra\E\ra\E'\ra O_D\ra0
\end{equation*}

An extension in the same orbit comes from $\Aut(O_D)$. As in the
proof of Lemma \ref{not_lolally_free}, let $V_i=\Ext(O_{p_i},\E)$.
Denote $L_i$ the line in $\PP(V)$ given by $\PP(V_i)$.

A point $p\in\PP(V)$ that does not lie on any of the lines
corresponds to a locally free extension. Note that the point $p$
is contained in a unique $(e-1)$-linear subspace in $\PP(V)$ that
intersects each line $L_i$ in only one point, say $q_i$. Since by
the action of $\Aut(O_D)$ we do not change the points $q_i$, it
follows that the orbit of the point $p$ is precisely the
$(e-1)$-linear subspace in $\PP(V)$ determined by the points
$q_i$. \ep


\subsection{The locus of unstable extensions}

\begin{lemma}\label{unstable_locus_2_local}
For any $D\in\Sym^e(C)$ and $\E\in M^s_{\xi(-D)}$ general, there is a closed
subvariety $Y_{D,\E}\subset\PP(V_{D,\E})\setm\Ga_{D,E}$, corresponding to
unstable extensions, and it has all the irreducible components of codimension
at least $2$. If $e=1$, then $Y_{D,\E}=\eset$.
\end{lemma}

\bp
The bundle $\E'$ in $(**)$ is not stable if and only if there is
a line bundle $\Lb$ on $C$ of degree $\deg(\Lb)\geq1$ and an
extension:
\begin{equation*}
0\ra\Lb\ra\E'\ra\Lb^{-1}\otimes\xi\ra0
\end{equation*}

Let $h$ be the morphism obtained by the composition $\Lb\ra\E'\ra O_D$.
If $h=0$, then $\Lb$ is a subbundle of $\E$. But $\E$ is a stable bundle
of degree $1-e$ and we get a contradiction, so we must have $h\neq0$. Then
there is $i\in\{1,\ldots,e\}$ and $D_i\in\Sym^i(C)$, with $D_i\subset D$,
such that $\im(f)=O_{D_i}$. Let $\K=\ker(h)$. It follows that there
are exact sequences:
\begin{equation}\tag{$**'$}
0\ra\K\ra\Lb\ra O_{D_i}\ra0
\end{equation}
\begin{equation*}
0\ra\K\ra\E\ra\K^{-1}\otimes\xi(-D)\ra0
\end{equation*}

Let $f=-\deg(\K)$. Note that we have $f=i-\deg(\Lb)\leq(e-1)$.
Since $\E$ is stable,  we have $f\geq\frac{e}{2}$.
Moreover, we have that $(**')$ maps to $(**)$  by the composition:
\begin{equation}\label{Ext_data}
\Ext^1(O_{D_i},\K)\ra\Ext^1(O_D,\K)\ra\Ext^1(O_D,\E)
\end{equation}

We have so far proved that when $\E$ is stable, we have that $\E'$
is unstable if and only if there are integers $i$ and $f$, such
that $e\geq i\geq1$ and $e-1\geq f\geq\frac{e}{2}$, a $0$-cycle
$D_i\subset D$ of length $i$, plus the following data: \bi
\item[i. ] a saturated line subbundle $\K$ of $\E$, with
$\deg(\K)=-f$ \item[ii. ] an element $(**')$ in
$\Ext^1(O_{D_i},\K)$, such that $(**)'$ is sent to $(**)$ via the
map (\ref{Ext_data}) \ei

We construct a parameter space for such data. Fix $f$ and $i$
integers such that $e-1\geq f\geq\frac{e}{2}$ and $e\geq i\geq1$
and fix $D_i\subset D$ a length $i$ $0$-cycle. Let
$q:X\ra\Pic^{-f}(C)$ be a projective bundle as defined in
(\ref{X}) and consider the morphism defined in (\ref{ka}):
\begin{equation*}
\ka:X\setm Z\ra M_{\xi(-D)}
\end{equation*}

If $\A$ is a Poincar\'e bundle on $\Pic^{-f}(C)\times C$ and
$\mu_1$, $\mu_2$ are the two projections from $X\times C$, then
consider the universal extension (\ref{univ_ext_1}) on $X\times C$:
\begin{equation*}
0\ra \mu_1^*O_X(1)\otimes q^*\A\ra\G\ra
q^*(\A)^{-1}\otimes\mu_2^*\xi(-D)\ra0.
\end{equation*}

Define $T$ to be an irreducible component of the fiber of $\ka$ at $\E$.
Consider the following relative extension sheaves on $T$:
\begin{equation*}
\F'=\EExt^1_{T\times C|T}(\mu_2^*O_{D_i}, q^*\A),\quad
\F=\EExt^1_{T\times C|T}(\mu_2^*O_{D}, \G)
\end{equation*}

The sheaves $\F'$ and $\F$ are locally free and
$\PP(\F')\subset\PP(\F)$ is a closed subvariety. Note that since
for any $x\in T$ we have $\G_x\cong\E$ and $\E$ is stable, it
follows that $\G\cong\mu_2^*\E\otimes\N$, for $\N$ some line
bundle from $T$.

Denote $V=V_{D,E}$, $Y=Y_{D,E}$. We have that $\PP(\F)\cong
T\times\PP(V)$. The unstable locus $Y$ is given by the union (over
all $f$, $i$, $D_i$ and irreducible components $T$) of the images
of $\PP(\F')$ in $\PP(V)$ via the map
\begin{equation*}
\PP(\F')\subset\PP(\F)\cong T\times\PP(V)\ra\PP(V)
\end{equation*}
It follows that all the irreducible components of $Y$ have dimension
at most
\begin{equation*}
\dim\PP(\F')=\dim T+\dim\Ext^1(O_{D_i},\K)-1=\dim T+i-1
\end{equation*}
If $\E$ is general, then we have that either $T$ is empty, or that
$T$ has dimension at most
\begin{equation*}
\dim X-\dim M_{\xi(-D)}=2f-e-g+2
\end{equation*}
It follows that all the irreducible components of $Y$ have
codimension in $\PP(V)$ at least
\begin{equation*}
\dim \PP(V)-(2f-e-g+i+1)=2(e-f)+(e-i)+g-2\geq2
\end{equation*}
\ep

\ps

For every integer $f$ such that $(e-1)\geq f\geq\frac{e}{2}$,
consider the morphism defined in (\ref{ka_global}):
\begin{equation*}
\tilde{\ka}:\X\setm\Z\ra M_d
\end{equation*}

If $\tilde{\ka}$ is dominant, let $M^0_{1-e}\subset M^s_{1-e}$ be the open
set where the fibers of $\tilde{\ka}$ have dimension $\dim\X-\dim M_{1-e}$.
If  $\tilde{\ka}$ is not dominant, we let $M^0_{1-e}$ be a dense open
in the complement of the image of $\tilde{\ka}$.
If we do this for all $f$, we end up with a dense open
$M^0_{1-e}\subset M^s_{1-e}$, so that if $\E\in M^0_{1-e}$,
then for any $D\in\Sym^e(C)$, $\E$ will be ``general'' in the sense
of Lemma \ref{unstable_locus_2_local}.

\ps

From now on, instead of $S'$ and $W'$ from (\ref{S'}) and
(\ref{W'}), we will work with
\begin{equation*}
S=S'\times_{M_{1-e}}M^0_{1-e},\quad W=W'\times_{S'}S
\end{equation*}

Let $Y\subset W\setm\Ga$ be the closed subset of $W$, corresponding to
those $w\in W$, for which the bundle $\J_w$ of (\ref{univ_ext_2}) is not
stable. Then we have the following Corollary to
Lemma \ref{unstable_locus_2_local}:
\begin{cor}\label{unstable_locus_2}
The subvariety $Y\subset W\setm\Ga$, corresponding to unstable extensions
in $W\setm\Ga$, has all irreducible components of codimension at least $2$.
If $e=1$, then $Y=\eset$.
\end{cor}

Denote $Z=Y\cup(\Ga\cap W)$ and if $Y_{D,\E}$ is as in Lemma
\ref{unstable_locus_2_local}, we let
\begin{equation*}
Z_{D,\E}=\Ga_{D,\E}\cup Y_{D,\E})
\end{equation*}

It follows from Corollary \ref{unstable_locus_2}, that the universal bundle
$\J$ of (\ref{univ_ext_2}), when restricted to $W\setm Z$,
induces a morphism:
\begin{equation}\label{eta}
\eta:W\setm Z\ra M
\end{equation}
which sends an extension $(**)$ to the isomorphism class of the bundle $\E$.



\section{The nice component of the space of rational curves of even degree}

Consider the morphism $\eta:W\setm Z\ra M$ of (\ref{eta}). Denote by
$\eta_{D,\E}$ the restriction of $\eta$ to the fiber of $p:W\ra S$:
\begin{equation*}
\eta_{D,\E}:\PP(V_{D,\E})\setm Z_{D,\E}\ra M
\end{equation*}

\ps

Note that by Lemma \ref{unstable_locus_2_local}, the codimension
of $Z_{D,\E}$ in $\PP(V_{D,\E})$ is at least $2$.



\subsection{The degree of the line bundle $\eta_{D,\E}^*\Th$}

\begin{lemma}\label{deg_2}
We have $\eta_{D,\E}^*\Th\cong O(2e)$.
\end{lemma}

\bp
Consider the restriction of the universal sequence (\ref{univ_ext_2})
to $\PP(V_{D,\E})\times C$:
\begin{equation*}
0\ra\mu_2^*\E(1)\ra\J_{D,\E}\ra\mu_2^*O_D\ra0
\end{equation*}
where $\J_{D,\E}$ is the restriction of the bundle $\J$ to
$\PP(V_{D,\E})\times C$. If we let $\{H\}$ be the class of a hyperplane
in $\PP(V_{D,\E})$, we have:
\begin{equation*}
c_1(\J_{D,\E})=2\{H\}\times C + \PP(V)\times c_1(\xi),\quad
c_2(\J_{D,\E})=\{H\}\times(c_1(\xi)+D)
\end{equation*}
The result follows by applying Lemma \ref{degree}.
\ep


\subsection{The nice component $\M$ of $\Hom_k(\PP^1, M)$ for $k=2e$}

\begin{lemma}\label{unobstr_2}
Let $f:\PP^1\ra M$ be the composition
\begin{equation}\label{comp_2}
\begin{CD}
\PP^1 @>h>>\PP(V_{D,\E})\setm Z_{D,\E} @>\eta_{D,\E}>> M
\end{CD}
\end{equation}
where $h$ is such that $h^*O(1)\cong O(1)$.
Then $f$ is an unobstructed point of $\Hom_k(\PP^1, M)$.
\end{lemma}

\bp
Note that by Lemma \ref{deg_2}, we have that $f^*\Th=O(k)$, so $f$ is
a point in $\Hom_k(\PP^1, M)$.
Consider the exact sequence obtained by pulling back the
universal extension (\ref{univ_ext_2}) to $\PP^1\times C$ by
$(h\times\id_C)$:
\begin{equation}\label{seq_2}
0\ra O(1)\bo\E\ra\F\ra O\bo O_D\ra0.
\end{equation}
The bundle $\F$ induces the morphism $f:\PP^1\ra M$.
It follows that if $\eta$ is the generic point of $C$, then
$\F_{\eta}\cong O(1)\oplus O(1)$. From the Lemma \ref{useful_result},
it follows that $f$ is an unobstructed point in $\Hom_k(\PP^1,M)$.
\ep

\begin{thm}\label{nice_comp_even}
If $k=2e$ is an even positive integer, there is an irreducible component
$\M$ of $\Hom_k(\PP^1,M)$, having the expected dimension $2k+3g-3$, with
the following properties:
\bi
\item[i)] A general point $f\in\M$ is unobstructed and it is
obtained as a composition:
\begin{equation*}
\begin{CD}
\PP^1 @>h>>\PP(V_{D,\E})\setm Z_{D,\E} @>\eta_{D,\E}>> M
\end{CD}
\end{equation*}
where $D\in\Sym^e(C)$, $\E\in M^s_{1-e}$ and $h$ is a morphism such
that $h^*O(1)\cong O(1)$
\item[ii)] The MRC fibration of $\M$ is given by a rational map:
\begin{equation*}
\M\dra\Pic^{e}(C)
\end{equation*}
which sends a point $f\in\M$ as in $i)$ to $O(D)\in\Pic^e(C)$. The
map is dominant if and only if $e\geq g$. \ei
\end{thm}

\bp
We follow the idea of the proof of Theorem \ref{nice_comp_odd}.

Let $\Hom_1(\PP^1,W\setm Z)$ be the relative Hilbert scheme morphisms
$\PP^1\ra W\setm Z$ that are contained in some fiber of $p:W\ra S$ and
have degree $1$ with respect to the relatively ample line bundle $O_W(1)$.
There is a canonical morphism:
\begin{equation}\label{si}
\si:\Hom_1(\PP^1,W\setm Z)\ra S
\end{equation}
whose fibers are the smooth, irreducible varieties
$\Hom_1(\PP^1,\PP(V_{D,\E})\setm Z_{D,\E})$. They have dimension
$(4e-1)$. If we let $N=\dim\PGL(r)+e-1$, we have:
\begin{equation*}
\dim\Hom_1(\PP^1,W)=\dim S+(4e-1)=(2k+3g-3)+N
\end{equation*}

We have that $\eta:W\setm Z\ra M$ induces a morphism
\begin{equation}\label{phi}
\phi:\Hom_1(\PP^1,W\setm Z)\ra\Hom_k(\PP^1,M)
\end{equation}

\ps

We prove that the general fiber of $\phi$ is $N$-dimensional and rational.
Let $h$ and $h'$ be two morphisms $\PP^1\ra W\setm Z$ such that
$\eta\circ h=\eta\circ h'$. Assume that $h$, respectively $h'$, have image
in some $\PP(V_{D,\E})$, respectively in $\PP(V_{D',\E'})$.
We have that the following compositions are equal:
\begin{equation*}
\begin{CD}
\PP^1 @>h>>\PP(V_{D,\E})\setm Z_{D,\E} @>\eta_{D,\E}>> M\quad\hbox{and}
\quad\PP^1 @>h'>>\PP(V_{D,\E})\setm Z_{D,\E} @> \eta_{D,\E}>> M
\end{CD}
\end{equation*}
As in (\ref{seq_2}) we have exact sequences on $\PP^1\times C$:
\begin{equation}\label{seq_2_F}
0\ra O(1)\bo\E\ra\F\ra O\bo O_D\ra0.
\end{equation}

\begin{equation}\label{seq_2_F'}
0\ra O(1)\bo\E'\ra\F'\ra O\bo
O_{D'}\ra0.
\end{equation}

Since $\F$ and $\F'$ induce the same morphism $\PP^1\ra M$, there
is an integer $m$ such that $\F'\cong\F\otimes p_1^*O(m)$ (without
loss of generality, we may assume $m\geq0$). It follows that there
is a commutative diagram:
\begin{equation}\label{diagram_2_F_F'}
\begin{CD}
0@>>>O(1)\bo\E@>>>\F@>>>O\bo O_D@>>> 0\\
@VVV @VVV  @VV{\cong}V @VVV @VVV\\
0@>>>O(1-m)\bo\E'@>>>\F'(-m)@>>>O(-m)\bo O_{D'}@>>>0
\end{CD}
\end{equation}

Since a surjective  morphism $O_D\ra O_D$ is an isomorphism, it
follows that all the vertical arrows in (\ref{diagram_2_F_F'}) are
isomorphisms. Hence, $m=0$, $\E\cong\E'$ and $\F\cong\F'$. But the
automorphisms of a stable bundle $\E$ are given by multiplication
by scalars; it follows that the two extensions (\ref{seq_2_F}) and
(\ref{seq_2_F'}) are in the same orbit of the action of
$\Aut(O_D)$ on $V_{D,\E}$. Moreover, the maps $h$ and $h'$ are in
the same orbit of the induced action of $\Aut(O_D)$ on
$\Hom_1(\PP^1,\PP(V_{D,\E}))$. It follows that the fiber of $\psi$
containing $h$ has dimension $N=\dim\PGL(r)+e-1$ and it is a
rational variety.

\ps

Let $\M$ be the closure of the image of $\phi$, with the
reduced structure. We have that $\M$ is an irreducible closed
subvariety of $\Hom_k(\PP^1,M)$, which has the expected
dimension $2k+3g-3$. By Lemma \ref{unobstr_2}, a general point $f\in\M$ is
unobstructed. It follows that $\M$ is the unique irreducible component of
$\Hom_k(\PP^1,M)$ containing the point $f$. This proves i).

\ps

For ii), consider the morphism $\si:\Hom_1(\PP^1,W\setm Z)\ra S$.
Let $\pi:S\ra\Pic^e(C)$ be the canonical morphism and let
$\th=\pi\circ\si$. We claim that $\th$ gives the MRC fibration
\begin{equation*}
\th:\Hom_1(\PP^1,W\setm Z)\ra\Pic^e(C)
\end{equation*}

It is enough to prove that the fibers of $\th$ are rationally connected.
Since the fibers of $\si$ are rational, using \cite{GHS}, it is enough to
prove that the fibers of $\pi$ are rationally connected.

Let $D\in\Sym^e(C)$. Then the fiber of $\pi$ at the point
$O_C(D)\in\Pic^e(C)$ is isomorphic to
\begin{equation*}
\tau^{-1}(M^s_{\xi(-D)})\times\PP(H^0(C,O(D))
\end{equation*}
The moduli spaces $M_{\xi(-D)}$ are unirational. As the fibers of
$\tau$ over the stable locus are isomorphic to $\PGL(r)$, it
follows that the fibers of $\pi$ are rationally connected.

Consider the morphism $\phi:\Hom_1(\PP^1,W\setm Z)\ra\M$. Since
the general fiber of $\phi$ is rationally connected, it follows that
$\th$ induces a morphism which gives the MRC fibration of $\M$:
\begin{equation*}
\rho:\M\dra\Pic^e(C)
\end{equation*}
\ep


\subsection*{Remark}
If $n\geq1$ is an integer, we let $k=2en$. Using the same arguments as
in the proof of Theorem \ref{nice_comp_odd}, one can still prove that there
is a morphism:
\begin{equation*}
\phi:\Hom_n(\PP^1,W\setm Z)\ra\Hom_k(\PP^1,M)
\end{equation*}

It follows that there are irreducible subvarieties
$\N(e,n)\subset\Hom_k(\PP^1,M)$, such that a general point $f\in\N(e,n)$
is obtained as a composition
\begin{equation}\label{comp_2'}
\begin{CD}
\PP^1@>h>>\PP(V_{D,\E})\setm Z_{D,\E}@>{\eta_{D,\E}}>>M
\end{CD}
\end{equation}
where $h$ is a morphism such that $h^*O(1)\cong O(n)$.

As in the proof of Theorem \ref{nice_comp_odd}, one can compute
the dimensions of the subvarieties $\N(e,n)$ and find that, if
$n>1$, these dimensions are strictly smaller than the expected
dimension $2k+3g-3$, so they cannot form irreducible components.
Note that for $n=1$ the subvariety $\N(k,1)$ is the nice component
$\M$.


\subsection{Geometric interpretation}

\

\ps

Let $D\in\Sym^e(C)$ and $\E\in M_{\xi(-D)}$ and  
denote $V=V_{D,\E}$, $Z=Z_{D,\E}$.
Assume that $D$ consists of distinct points $p_1,\ldots,p_e$. We let
$V_i=\Ext^1(O_{p_i},\E)$ and consider the lines $L_i=\PP(V_i)$ in
$\PP(V)$. Let $l$ be a line in $\PP(V)\setm Z$. 
By projecting from the linear span of all the lines other than $L_i$, we 
obtain isomorphisms $l\ra L_i$, that determine
a rational normal scroll $\Si\subset\PP(V)$. Then 
$\Si\cong\PP^1\times\PP^{e-1}$ and the fibers of $\Si$ over
$\PP^1$ intersect the lines $L_i$ exactly once and
$l$ is a section of $\Si\ra\PP^1$. 

\ps

Consider the action of $\Aut(O_D)$ on the space of lines in $\PP(V)\setm Z$.  
Then the orbit containing the line $l$ is given by the sections of 
$\Si\ra\PP^1$. Note that they form an $(e-1)$-dimensional family.


\section{Free rational curves of even degree}

\begin{lemma}\label{free_2}
A morphism $f$ which is a composition as in (\ref{comp_2}) is a free curve.
\end{lemma}

\bp
We follow the same arguments as in Lemma \ref{free_1}. The morphism $f$
is induced by a vector bundle $\F$ on $\PP^1\times C$, which sits in an
exact sequence:
\begin{equation}\label{seq_F_2}
0\ra O(1)\bo\E\ra\F\ra O\bo O_D\ra0
\end{equation}

Let $\T=\R^1{p_1}_*(\F\otimes\F^*)$ and recall that $\T\cong O^g\oplus f^*T_M$.
We prove that $\T$ is non-negative. 

\ps

Denote
\begin{equation*}
\S'=\F\otimes(O(-1)\bo\E^*)\quad\quad\S''=\F\otimes(O\bo O_D)
\end{equation*}

\ps

By tensoring (\ref{seq_F_2}) with $O(-1)\bo\E^*$, respectively with
$O\bo O_D$, we get exact sequences:

\begin{equation}\label{1}
0\ra O\bo(\E\otimes\E^*)\ra\S'\ra O(-1)\bo(\E^*\otimes O_D)\ra0
\end{equation}

\begin{equation}\label{2}
0\ra O\bo O_D\ra O(1)\bo(\E\otimes O_D)\ra\S''\ra O\bo O_D\ra0
\end{equation}

Note that from (\ref{2}) it follows that 
\begin{equation}\label{S''}
\S''\cong (O(2)\oplus O)\bo O_D
\end{equation}

This implies that ${p_1}_*(\S'')\cong O(2)^e\oplus O^e$ and 
$\R^1{p_1}_*(\S'')=0$.

\ps

Note that $h^0(C,\E^*\otimes O_D)=2e$ and $h^1(C,\E^*\otimes O_D)=0$.
The long exact sequence coming from applying ${p_1}_*$ to
(\ref{1}) has the form:
\begin{equation}\label{1'}
\begin{CD}
0\ra O\otimes\H^0(C,\E\otimes\E^*)\ra{p_1}_*(\S')\ra O(-1)^{2e}\ra \\
\ra O\otimes\H^1(C,\E\otimes\E^*)\ra\R^1{p_1}_*(\S')\ra0
\end{CD}
\end{equation}

It follows from (\ref{1'}) that the sheaves ${p_1}_*(\S')$ and 
$\R^1{p_1}_*(\S')$ are locally free and that $\R^1{p_1}_*(\S')$ is 
non-negative.

\ps

By dualizing (\ref{seq_F_2}) and tensoring with $\F$, we get an exact 
sequence:
\begin{equation}\label{3}
0\ra \F\otimes\F^*\ra\S'\ra\S''\ra0
\end{equation}

Applying ${p_1}_*$ to (\ref{3}) and using that 
${p_1}_*(\F\otimes\F^*)\cong O$ (and  (\ref{S''})), we have a sequence

\begin{equation}\label{3'}
0\ra O\ra{p_1}_*(\S')\ra O(2)^e\oplus O^e\ra\T\ra\R^1{p_1}_*(\S')\ra0
\end{equation}
As $\R^1{p_1}_*(\S')$ is non-negative, it follows from (\ref{3'})
that $\T$ is non-negative.
\ep

\begin{cor}
Let $k$ be an even integer and let $\M\subset\Hom_k(\PP^1,M)$ be the nice
component. If $k$ is sufficiently large, a general $f\in\M$ is a very 
free curve.
\end{cor}

\bp
We make a similar argument as in the odd degree case. Let $k=2e$ and let
$D\in\Sym^e(C)$, $\E\in M_{\xi(-D)}$.
Using Lemma \ref{dominant}, if $e$ is sufficiently large, the morphism 
\begin{equation*}
\eta_{D,\E}:\PP(V_{D,\E})\setm Z_{D,\E}\ra M
\end{equation*}
is dominant for general $\E$ and $D$. 
The result now follows from Lemma \ref{very_free}, using the same argument
as in the odd degree case, see Corollary \ref{very_free_odd} iii). 
\ep

\begin{lemma}\label{dominant}
If $e$ is sufficiently large, the morphism 
\begin{equation*}
\eta_{D,\E}:\PP(V_{D,\E})\setm Z_{D,\E}\ra M
\end{equation*}
is dominant, for general $D\in\Sym^e(C)$ and  $\E\in M_{\xi(-D)}$.
\end{lemma}

\bp
Recall that if $\ze$ is a line bundle on $C$, we denote by $M_{\ze}$
the moduli space of semistable, rank $2$ vector bundles on $C$, with 
determinant $\ze$.

\ps

To prove the lemma, note that it is enough to prove the following:
if $\xi$ and $\ze$ are line bundles of degrees $1$, respectively
$(1-e)$, for $e$ sufficiently large, the determinant map
\begin{equation}\label{det}
\Hom(\E,\E')\ra\Hom(\wedge^2(\E),\wedge^2(\E'))\cong\H^0(\xi\otimes\ze^{-1})
\end{equation}
is dominant for $(\E,\E')$ general in $M_{\ze}\times M_{\xi}$. This is 
because if $h:\E\ra\E'$ is an injective morphism, 
that maps by the determinant map (\ref{det}) to a section in  
$\H^0(\xi\otimes\ze^{-1})$ whose divisor of zeros is an effective 
divisor $D\in\Sym^e(C)$, then we have $\coker(h)\cong O_D$.
Hence, $\E'$ is in the image of the map $\eta_{D,\E}$.
\ps

We prove the statement by showing that if $\E$ and $\E'$ are some 
well-chosen direct sums of line bundles on $C$, the determinant map is 
surjective. By deforming $\E$ and $\E'$ to stable vector bundles, the general
statement follows. Note that if $e$ is sufficiently large, we have
\begin{equation}\label{h^0}
\h^0(\E^*\otimes\E')=\chi(\E^*\otimes\E')=2e+4-4g
\end{equation}
for any $\E\in M_{\ze}$, $\E'\in M_{\xi}$.
This is because  there are no nonzero morphisms $\E'\ra\E\otimes K_C$
if $e$ is sufficiently large. Otherwise, such a morphism should factor 
through a line bundle, which would then lead to a contradiction of the 
stability of $\E$ and $\E'$. It follows that
\begin{equation*}
\H^0(\E\otimes{\E'}^*\otimes K_C)\cong\Hom(\E',\E\otimes K_C)=0
\end{equation*}
which then implies (\ref{h^0}).

\ps

Let $p\in C$ and let $\eta\in\Pic^0(C)$ such that $\eta^2\cong\xi(-p)$.
Let $d$ be an integer such that $(e-1)\geq d\geq(\frac{e}{2}+g)$ and let
$\Lb$ be a line bundle on $C$ of degree $d$. Let 
$\M=\ze^{-1}\otimes\xi\otimes\Lb^{-1}(-p)$ and we consider the following
rank $2$ vector bundles:
\begin{equation*}
\E=(\Lb^{-1}\otimes\eta)\oplus(\M^{-1}\otimes\eta), \quad\quad
\E'=\eta\oplus\eta(p)
\end{equation*}

There is an isomorphism
\begin{equation*}
\H^0(\E^*\otimes\E')\cong
\H^0(\Lb)\oplus\H^0(\M)\oplus\H^0(\Lb(p))\oplus\H^0(\M(p))
\end{equation*}

\ps

The determinant map sends an element $(s,t,s',t')$ to the element
$(st'-s't)\in\H^0(\Lb\otimes\M(p))$, where 
$s\in\H^0(\Lb)$, $t\in\H^0(\M)$, $s'\in\H^0(\Lb(p))$, $t'\in\H^0(\M(p))$.
We prove this map is surjective.

\ps

By the base-point-free pencil trick (see \cite{ACGH}, III.3), 
if we fix $t'_1$ and $t'_2$ sections in $\H^0(\M(p))$ that do not have a 
common zero, then any section in $\H^0(\Lb\otimes\M(p))$ can be written as 
$t'_1s_1-t'_2s_2$, for some $s_1, s_2$ sections in $\H^0(\Lb)$, provided
that $H^1(\Lb\otimes\M^{-1}(-p))=0$. Note that our choice of $d$ was exactly
so that we have this vanishing. Note that we may choose
$t'_2$ to be of the form $rt$, where $r$ is a non-zero element in 
$\H^0(O(p))$ and $t\in\H^0(\M)$. Then 
the element $(s_1, t, rs_2,t'_1)$ will map by the determinant map to 
$t'_1s_1-t'_2s_2$. We have proved that the determinant map is surjective
and the lemma follows.
\ep

Note that by Lemma \ref{fibers}, the fibers of the morphism $\eta_{D,\E}$ 
have dimension at least $(e-1)$. It follows that the condition
$e\geq(3g-3)$ is necessary for the morphism $\eta_{D,\E}$ to be dominant. 
At best, the previous method would only give the bound $k\geq(6g-6)$ for the 
curves on $M$ in the nice component (of even degree) to be very free. 
We may improve this bound by producing very free rational curves of even 
degree as deformations of reducible curves which are unions of two
free curves of odd degree, intersecting at a point, of which at least one 
is very free. (Using a fact that we prove in \cite{C}, namely that 
all the irreducible components other than the nice component and the 
almost nice component are obstructed, it follows that any very free curves 
will have to be in the nice component.)


\section{Example -- the genus $2$ case}

If $g=2$, recall that $M$ is isomorphic to a complete 
intersection of two quadrics in $\PP^5$.

\ps

For any $k\geq1$, the space $\Hom_k(\PP^1,M)$ has two irreducible 
components of the expected dimension $2k+3$:
\bi
\item[i) ]the nice component $\M$, whose general element is a
is very free if $k$ is sufficiently large
\item[ii) ]the almost nice component $\M'$, whose general element is 
a morphism $f:\PP^1\ra M$ which is $k$-to-$1$ onto a line in $M$
\ei
(The two components are the same if $k=1$.)
The MRC fibration of both components is given by a rational map to $J(C)$,
which is dominant for all $k$, except $k=2$.


\end{document}